\documentclass[12pt]{article}
\usepackage[applemac]{inputenc}              % afficher le caractère français sur mac
\usepackage{amsfonts}
\usepackage{amsmath}
\usepackage{graphicx}
\usepackage[francais]{babel}
\usepackage{color}
\newtheorem{remark}{Remarque}[section]
\newtheorem{definition}{D\'efinition}[section]
\newtheorem{theorem}{Th\'eor\`eme}[section]

\newtheorem{corollary}{Corollaire}[section]
\renewcommand{\sectionmark}[1]{\markright{\underline{ I. Faye, I. Lavall\'ee, M. Ngom,  D. Seck, A. Sy}.}}

\title{{\LARGE Programmation Lin\'eaire, une nouvelle approche}\\ {\Large Novel way in linear programming}}
\author{ I. Faye$^{1,2}$\, ;  I. Lavall\'ee$^{4}$\footnote{ivan.lavallee@gmail.com}\, ; M. Ngom$^{1,2}$\, ; D. Seck$^{2,3}$\, ; A. Sy$^{1,2}$\footnote{\{azousy2, dseck, grandmbodj, ngomata\}@hotmail.com}}

\begin{document}
\maketitle
  {\small \centerline{$^{1}$ U.F.R de Sciences Appliqu\'ees et des Technologies de}
   \centerline{l'Information et de La Communication} 
   \centerline{Universit\'e de Bambey BP 30, Bambey, S\'en\'egal}
%   \medskip
   \centerline{$^{2}$ Laboratoire de Math\'ematiques de la D\'ecision et d'Analyse Num\'erique}
 \centerline{\'Ecole Doctorale Math\'ematiques et Informatique}
   \centerline{UCAD BP 16 889 Dakar-Fann S\'en\'egal}
 \centerline{$^{3}$ Facult\'e des Sciences Economiques et de Gestion  UCAD. Dakar}}
 \medskip

\centerline{$^{4}$ LaISC Université Paris 8 \& C.N.R.S. UMI ESS 3189, BP. 5005, UCAD Dakar}
\begin{abstract}
Après un bref aperçu permettant de situer notre travail, nous proposons une nouvelle voie pour aborder la programmation linéaire en proposant un algorithme élaboré à partir d'une idée simple qui permet d'obtenir une solution aussi approchée que voulu par translation dichotomique d'un hyperplan de l'espace des solutions.
\end{abstract}
\textbf{Mots clés:}\quad {\small Algorithme, Programmation Linéaire, polyèdre, polytope, hyperplan, dichotomie, optimisation, approximation.}
\begin{center}
\begin{minipage}{11.5cm}
\small {\center \textbf{Abstract} \par}
\quad After a short course in order to situate our work, we propose a new way to study linear programming and we give a proposal of algorithm to solve linear programming from a basic idea which allow to obtain an approached solution with desired accuracy. For this we use some dichotomic translations of an hyperplan in the solutions hyperspace.
\end{minipage}
\end{center}
\section{Introduction}
Le problème de programmation linéaire s'est rapidement imposé dès qu'on a voulu planifier un tant soit peu les activités économiques ou autres. C'est ainsi que dès les années 1939, les nécessités de la planification soviétique conduisent Kantorovitch et Tolstoï à proposer une solution au problème, (voir \cite{KAN39,TOL39}) inspirée plus ou moins des travaux de Joseph Fourier (1768-1830). En 1951 Dantzig publie les résultats de Kantorovitch-Tolstoï sous forme d'algorithme exécutable sur ordinateur (voir \cite{DAN51}) et lui donne le nom de \textit{simplex}. Cet algorithme fait le tour du monde sans aucun concurrent jusqu'en 1979 avec l'apparition de l'algorithme dit de \textit{l'ellipsoïde} (voir \cite{KAT79}), qui est en \textit{temps polynomial}, ce qui, dans le pire des cas (dégénéres\-cence), n'est pas le cas du \textit{simplex}. Paraît ensuite un autre algorithme polynomial basé sur une autre idée, celui de Karmarkar (voir \cite{KAR84-a} et \cite{KAR84-b}), puis en 1989 l'algorithme de Murty \& Chang (voir \cite{MUR89}). Ce dernier algorithme étant amélioré en 2006 (voir \cite{MUR06-a,MUR06-b}). Ces algorithmes, {\em simplex, ellipsoïde, Karmarkar, Murty} sont souvent caractérisés comme étant, {\em de frontière} pour le \textit{simplex}, \textit{extérieur} pour {\em l'algorithme de l'ellipsoïde}, \textit{intérieur} pour {\em les algorithmes de Karmarkar et de Murty}. D'autres aussi proposent une approche un peu différente (voir par exemple \cite{NPL09, PH09}). Dans le présent article, nous avons voulu présenter une méthode, qu'on peut qualifier, d'\textit{extérieure-intérieure}, ou encore \textit{alternante} basée sur une idée géométrique de dichotomie de l'espace. Il s'agit de considérer (dans le cas d'une maximisation) un hyperplan particulier de la famille de ceux définis par la fonctionnelle à optimiser, tantôt extérieur au polytope des contraintes, tantôt le recoupant, et de se rapprocher par translation contrôlée, autant qu'on le souhaite de la solution optimale par un argument de dichotomie. On obtient ainsi une solution approchée aussi précise que voulu, ce qui est essentiel dans le domaine pratique, l'optimum sur le modèle n'étant quasiment jamais l'exact optimum dans la réalité modélisée. 

\section{Position du probl\`eme}
Dans toute la suite E d\'esigne un $\mathbb{R} -$espace-vectoriel; $\mathbb{R}$ étant l'ensemble des nombres réels.

On appelle probl\`eme de \textit{programmation lin\'eaire}, tout probl\`eme pouvant s'\'ecrire sous la forme:

\begin{eqnarray}\label{proli}
  &\max \sum_{i=1}^{n} c_{i}x_{i} \quad c_i \in \mathbb{R}; n \in \mathbb{N}^* & \notag \\
  s/c &\left\{\begin{array}{ccc}
               \sum_{i=1}^{n} a_{ij}x_{j} & \leq & b_{i} \\
               x_{j} & \geq & 0
             \end{array}
  \right.\qquad a,b \in \mathbb{R}; \; j=1:m,i=1:n&
\end{eqnarray}

On rappelle qu'un poly\`edre convexe de $E$ est un ensemble $P$ de la forme
\begin{equation}\label{poly1}
    P=\{x\;\in \;\mathbb{R}^{n}\;:\; Ax\;\leq\;b\;\},
\end{equation}
o\`u $A\;:\; E\,\longrightarrow \mathbb{R}$ est une application lin\'eaire ($m\;\in \,\mathbb{N}$;
si $m=0,\;P=E$,) $b\;\in \mathbb{R}^{m}$ et l'in\'egalit\'e $Ax\leq b$ se lit composante par composante dans $\mathbb{R}^{m}$
$(Ax)_i\leq b_i$ pour tout $i\in\{1,\cdots, m\}.$ Si l'ensemble se pr\'esente avec des \'egalit\'es lin\'eaires,
$Cx=,d$ on pourra se ramener sous la forme (\ref{poly1}) en les rempl\c{c}ant par deux in\'egalit\'es oppos\'ees $Cx\leq d$ et
$-Cx\leq -d.$ Un polytope est un poly\`edre convexe et born\'e.

G\'eom\'etriquement, un poly\`edre est donc l'intersection d'un nombre fini de demi-espaces de $E$. Si $E$ est de dimension finie, il n'y a
pas de restriction \`a supposer que $E=\mathbb{R}^{n}$ et que $A$ est une matrice $m\times n$ (il suffit de donner une base de $E$).
Dans certaines circonstances, par exemple en optimisation lin\'eaire, il est plus avantageux de repr\'esenter un poly\`edre de $\mathbb{R}^{n}$
sous la forme dite standard suivante:
\begin{equation}\label{poly2}
    P=\{x\;\in\; \mathbb{R}^{n}\; Ax=b,\;x\geq 0\;\}.
\end{equation}
Il n'y a aucune perte de g\'en\'eralit\'e dans cette repr\'esentation. Tout poly\`edre de la forme (\ref{poly1}) se repr\'esente sous la forme
(\ref{poly2}) en introduisant des variables d'\'ecart.
\begin{remark}
Les repr\'esentations (\ref{poly1}) et (\ref{poly2}) sont dites duales, car elles font intervenir des applications lin\'eaires (\'el\'ements du
dual de $E$).
\end{remark}
On consid\`ere le polytope $D$ (voir figure) suppos\'e contenir l'origine $O$,
 sinon, on peut toujours introduire les variables d'\'ecart de sorte que l'origine $O\in D.$\\
%\begin{figure}
  \includegraphics[width=4in]{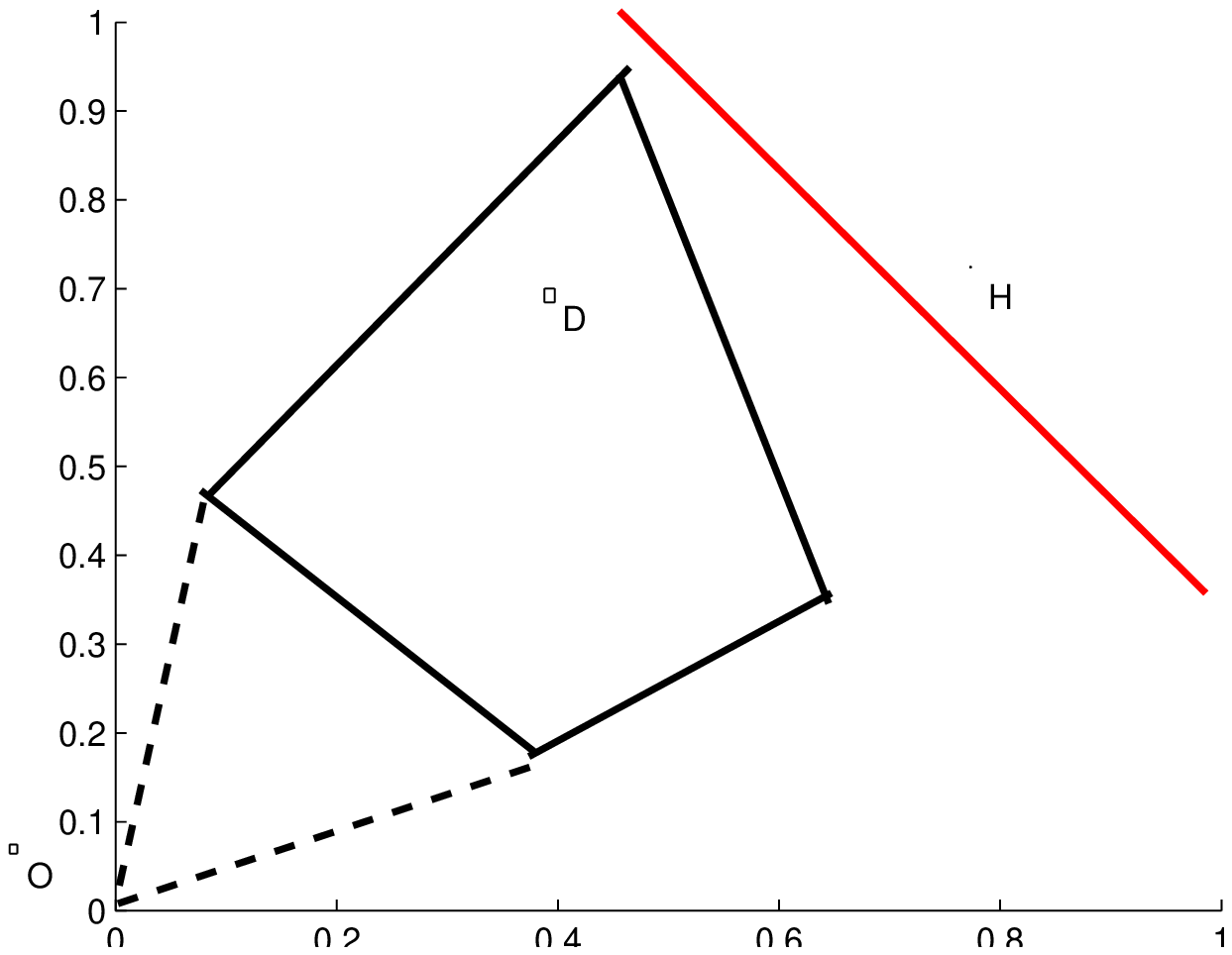}\\
  \textit{S\'eparation de l'espace en deux parties par un hyperplan dont l'une contient le polytope.}\label{separe}
%\end{figure}

\section{Quelques rappels en optimisation}
Soient  $f$ une forme lin\'{e}aire sur $E$, et $\alpha \in \mathbb{R} ,A,B$ deux sous
ensembles de $E.$
$$f(X)=\sum_{j=1}^{n}k_jx_j,\; n\in \mathbb{N^*}; k_j\in \mathbb{R}$$
\begin{definition}\label{def1}
 Un sous-ensemble C de E est dit convexe
si $\forall x,y\in C$ et $\lambda \in \lbrack 0,1]$ on a $\lambda
x+(1-\lambda )y\in C$ , On note $$[x,y]=\{\lambda x+(1-\lambda )y,0\leq
\lambda \leq 1\}$$
\end{definition}

\begin{definition}\label{def2} On appelle hyperplan $H$ d'\'{e}quation $[f=\alpha ]$, l'ensemble
d\'efinit par
\begin{equation}\label{hyplan}
    H=\left\{X\;\in\;\mathbb{R}^{N},\;\; tel\;\;que\;\; f(X)=\alpha. \right\}
\end{equation}

\end{definition}
$H$ est ferm\'{e} $\Longleftrightarrow f$
est continue.
%%%%%%%%%%%%%%%%%%%%%%%%%%%%%%%%%%%%%%%%%%%%%%%%%%%%%%%%%%%%%%%%%%%%%%%%%%%%%%%%%%%%%%%%%%%%%%%%%%%%%%%%%%%%%%%%%%%%%%%%%%%%%%%%%%%
\subsection{S\'eparation des ensembles convexes}
Un outil essentiel en analyse convexe est le théor\'eme de Hann-Banach sur la s\'eparation des ensembles convexes.
On supposera que l'espace vectoriel $E$ est de dimension finie. On peut toujours le munir d'un produit scalaire not\'e
$\langle .,.\rangle.$

La s\'eparation de deux ensembles convexes se fait g\'eom\'etriquement dans $E$ en utilisant un hyperplan affine $H$, c'est-\`a-dire
un ensemble de la forme
$$H =\{ x\;\in\;E\;: \; \langle \xi, x\rangle = \alpha\;\};$$
o\`u $ \xi \in E$ est non nul et $\alpha \in \mathbb{R}$. On dit que cet hyperplan s\'epare deux convexes $C_1$ et $C_2$ si l'on a
$$\forall\; x_1\;\in\; C_1,\;\; \forall\; x_2\;\in\; C_2\;:\; \langle \xi, x_1\rangle \leq \alpha\leq \langle \xi, x_2\rangle $$
Ceci sera certainement le cas s'il existe un $\xi$ non nul dans $E$ tel que
$$\sup_{x_1\in C_1}\langle \xi, x_1\rangle \leq \inf_{x_2\in C_2}\langle \xi, x_2\rangle$$
On dit que cet hyperplan s\'epare strictement ces deux convexes s'il existe deux scalaires $\alpha_1$ et $\alpha_2$ tels que
$\alpha_1<\alpha < \alpha_2$ et
$$\forall\; x_1\;\in\; C_1,\;\; \forall\; x_2\;\in\; C_2\;:\; \langle \xi, x_1\rangle \leq \alpha_1<\alpha_2 \leq \langle \xi, x_2\rangle $$
Ceci sera certainement le cas s'il existe un $\xi$ (n\'ecessairement non nul) dans $E$ tel que
$$\sup_{x_1\in C_1}\langle \xi, x_1\rangle < \inf_{x_2\in C_2}\langle \xi, x_2\rangle$$

La figure ci dessous illustre cette notion. Ici, on a utilis\'e le produit euclidien sur $\mathbb{R}^{2}$.

\begin{center}
  \includegraphics[width=4in]{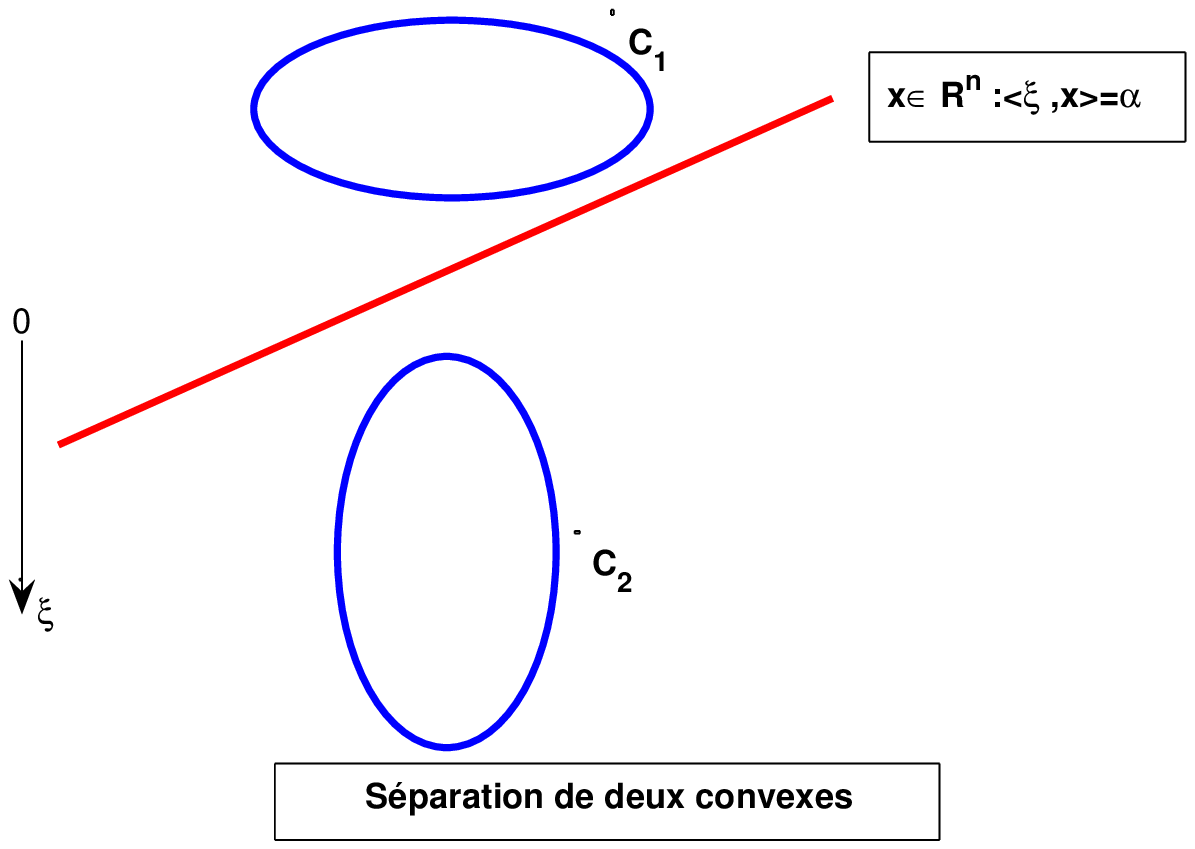}\\
  %\caption{}\label{}
\end{center}
\subsection{Th\'eor\`emes de Hann-Banach}
Nous allons donner ici quelques r\'esultats classiques en optimisation convexe et dont les d\'emonstrations peuvent \^etre trouv\'ees dans \cite{GIL06}, et aussi dans \cite{ACH84} pour ce qui concerne les polyèdres. Ces r\'esultats sont connus
 dans la litt\'erature sous le nom de formes g\'eom\'etriques du th\'eor\`eme de Hann-Banach
  ou de s\'eparation des ensembles convexes. Le premier affirme que l'on peut s\'eparer strictement deux ensembles convexes disjoints, si l'un est ferm\'e et l'autre est compact. Le second exprime que l'on peut s\'eparer (non strictement cette fois) deux ensembles convexes quelconques (en dimension finie).
\begin{theorem}[Hann-Banach I] Soit $E$ un espace vectoriel de dimension finie, muni d'un produit scalaire not\'e $\langle.,.\rangle$. Soient $C_1$ et $C_2$ deux convexes disjoints non vides de $E$, tels que $C_{1}^{\infty}\cap C_{2}^{\infty}=\{0\}$.
Alors, on peut s\'eparer $C_1$ et $C_2$ strictement: il existe un vecteur $\xi\;\in\;E$ tel que
$$\sup_{x_1\in C_1}\langle \xi, x_1\rangle < \inf_{x_2\in C_2}\langle \xi, x_2\rangle$$
\end{theorem}
\begin{corollary}
Soient $C_1$ et $C_2$ deux convexes non vides disjoints d'un espace vectoriel $E$ muni d'un produit scalaie not\'e
$\langle.,.\rangle$, on suppose l'un ferm\'e et l'autre compact, alors on peut s\'eparer $C_1$ et $C_2$ strictement.
\end{corollary}
\begin{theorem}[Hann-Banach II] Soit $E$ un espace vectoriel de dimension finie, muni d'un produit scalaire
not\'e $\langle.,.\rangle$. Soient $C_1$ et $C_2$ deux convexes disjoints non vides de $E$.
Alors, il existe un vecteur $\xi\;\in\;E$ non nul tel que
$$\sup_{x_1\in C_1}\langle \xi, x_1\rangle \leq \inf_{x_2\in C_2}\langle \xi, x_2\rangle$$
\end{theorem}
\begin{definition}
Soient $C\subset E$ un convexe et $H$ un hyperplan ferm\'{e} tel que $\exists u\in C\cap H,$ on
peut supposer que $H=\{f=\alpha \}$ $f\in E^{\prime },$ $E'$ \'etant le dual de $E$, et $\alpha \in \mathbb{R} .$ si $C\leq \{f\geq \alpha \}$ ou $%
C\subset \{f\leq \alpha \}$ on dit que $H$ supporte $A$ en $u$ ou que $u$ est un
point de support de $C$.
\end{definition}

\begin{corollary}
1- Soit $C$ un convexe ferm\'{e} non vide
$\forall u\in \partial C$ (bord de $C$), u est un point de support de $C$; Si
$C$ est un ensemble convexe d'un espace vectoriel topologique, par tout point
fronti\`{e}re passe un hyperplan d'appui.

2- Soit $C$ un ensemble convexe ferm\'{e}, Alors $$C=\bigcap\limits_{\Pi \in p}\Pi $$
$p$ est\;l'ensemble\;des\; demi-espaces\; ferm\'es\; contenant $C$.
\end{corollary}

\section{Algorithme propos\'e}
Dans ce paragraphe, nous proposons un algorithme d'optimisation permettant
de r\'esoudre le probl\`eme de programmation lin\'eaire donn\'e sous forme (\ref{proli}).

Soit $H$ l'hyperplan d\'efini par (\ref{hyplan}) et $K$ le polytope d\'efini  par
l'ensemble des contraintes de (\ref{proli}). En utilisant les th\'eor\`emes de
s\'eparation, on peut toujours s\'eparer l'espace en deux parties par l'hyperplan
 de sorte qu'une partie contienne le polytope $K$ des contraintes. 
 
 Trois possibilités apparaissent alors :
 \begin{itemize}
\item soit l'hyperplan $H$ ne coupe pas le polytope $K$ des contraintes et dans le cas d'une maximisation comme dans le PL de la formule \ref{proli} et alors $H$ est \og{} au-delà\fg{} des solutions réalisables, c'est-à-dire qu'il n'en peut contenir aucune et il convient alors de choisir un autre hyperplan $H^{'} // H$ situé entre l'origine du repère (\textit{i.e.} le point $0$) et $H$;
\item soit l'hyperplan $H$ coupe le polytope des contraintes et il faut alors inclure $H$ dans les contraintes et choisir un autre hyperplan de séparation de l'espace $H^* // H$ situé dans le demi-espace ne contenant pas $0$;
\item soit le polytope des contraintes est vide.
\end{itemize}
Si $H\cap K=\emptyset$, alors, en reprenant la notation de la définition \ref{def2}, on d\'efinit un nouvel hyperplan $H'=\{f=\frac{\alpha}{2}\}$
par exemple ; il va exister alors deux valeurs de $\alpha$, $\alpha_{s_0}$ et $\alpha_{s_1}$ , avec $s_{-1} \text{ et }s_0  \in
\mathbb{N}^{*}$  telles que l'hyperplan $H_{s_0}=\{f=\alpha_{s_0}\}\cap K\neq \emptyset$ alors que $H_{s_{-1}}=\{f=\alpha_{s_{-1}}\}\cap K = \emptyset$. On pose alors:
$$H_{s_{+1}}=\{f=s_{+1}=\dfrac{\alpha_{s_{0}}+\alpha_{s_{-1}}}{2}\}$$
Puis si $H_{s_{+1}} \bigcap K = \emptyset$, on pose $s_{-1}=s_{+1}$;

Sinon, si $H_{s_{+1}} \bigcap K \neq \emptyset$, on pose $s_{0}=s_{+1}.$

Dans cet algorithme, l'hyperplan choisi est d\'efini par
\begin{equation}\label{hyper3}
H=\{x\;\in\;\mathbb{R}^{n}\;/ \;\sum_{i=1}^{n}c_ix_i=\alpha\}
\end{equation}
o\`u $\sum_{i=1}^{n}c_ix_i$ repr\'esente la fonction objectif du probl\`eme de programmation lin\'eaire.
\begin{remark}
A chaque pas de l'algorithme, le nouvel hyperplan se d\'eduit du pr\'ec\'edent par translation.
En effet, on cherche les vecteurs qui engendrent l'espace $H$ d\'efini par (\ref{hyper3}) puis on cherche
un vecteur unitaire $\overrightarrow{v}$ qui repr\'esente le pas vectoriel de la translation......
\end{remark}
\subsection{Un sous-problème clé}
Dans cet algorithme, une difficulté réside dans l'identification de la vacuité -ou la \textbf{non}-vacuité- de l'intersection de l'hyperplan avec le polyèdre des contraintes :
\begin{equation}
H \bigcap K=\emptyset \, ?
\end{equation}
Au plan théorique, l'article de Murty \cite{MUR06-a} résoud le problème polynomialement. Du point de vue pratique, dans les cas de non-dégénérescence, le problème peut être résolu par l'utilisation de la méthode dite "\textit{du grand M}". 
%\section{Exemples Num\'eriques}

\section{Conclusion et Perspectives}\footnote{Cette étude nous a été facilitée par l'amitié et l'aide bibliographique de Pham Canh Duong du CNST à Hanoï.}
Cette étude ouvre une voie originale en matière de programmation linéaire qui n'est inspirée d'aucune des méthodes existantes mais qui utilise toutefois la méthode dit \textit{du point intérieur} de Murty pour la résolution d'un sous problème clé.
Il reste à examiner en quoi cette méthode est utilisable en {\emph programmation linéaire en nombres entiers} et en {\emph programmation linéaire paramétrique}. La méthode peut aussi s'étendre à la programmation non linéaire dans quelques cas particuliers.

%\begin{thebibliography}{99}
%\bibitem{opt} J. C. Gilbert \textit{OPTIMIZATION DIFF\'ERENTIELLE: Th\'eorie et Algorithme,} Cours \'Ecole Nationales Sup\'erieurede Techniques Avanc\'ees, 2005-2006.
%\newpage
%\vspace*{2cm}
%\hspace*{-3cm}{\rm \scalebox{4}[4]{$\phi = \frac{1}{2} \cdot (\sqrt{5} +1)= 1,618\ldots$} }
%$$ \phi = \frac{1}{2} \cdot (\sqrt{5} +1)= 1,618 \ldots$$
%\newpage
%{\rm \scalebox{2.5}[2.5]{$ \rho \left [ \partial_t \vee + (\vee.\nabla)\vee \right ]=- \nabla p + \nabla . \text{\bf T+ f}$}

%\end{thebibliography}
\bibliographystyle{alpha}
%\clearemptydoublepage
\bibliography{bib_ivan-tout}

\begin{thebibliography}{Mur06b}

\bibitem[Ach84]{ACH84}
S.~Achmanov.
\newblock {\em Programmation Linéaire}.
\newblock {\'E}dition Mir, 1984.
\newblock Traduit du Russe, édition Russe 1981.

\bibitem[CM89]{MUR89}
S.~Y. Chang and K.~G. Murty.
\newblock The steepest descent gravitational method for linear programming.
\newblock {\em Discrete Applied Mathematics}, 25:211--239, 1989.

\bibitem[Dan51]{DAN51}
G.B. Dantzig.
\newblock Minimization of a linear function of variables subject to linear
  inequalities.
\newblock In T.C. Koopman, editor, {\em In Actvity, Analysis of production and
  allocation}, pages 339--347. John Wiley, New York, 1951.

\bibitem[Gil06]{GIL06}
J.C. Gilbert.
\newblock {\em Optimisation diff{\'e}rentielle: Th{\'e}orie et algorithmes}.
\newblock {\'E}cole Nationale Sup{\'e}rieure des Techniques Avanc{\'e}es,
  2005-2006.
\newblock Cours.

\bibitem[Kan39]{KAN39}
L.V. Kantorovitch.
\newblock Méthodes mathématiques d'organisation et planification de la
  production, 1939.
\newblock En Russe Léningrad 1939, traduction anglaise Mathematical methods in
  the optimization and planning of production Management Science, Vol. 6,
  363-422, 1960.

\bibitem[Kar84a]{KAR84-a}
N.~Karmarkar.
\newblock A new polynomial-time algorithm for linear programming.
\newblock {\em STOC}, 1(11):302--311, 1984.

\bibitem[Kar84b]{KAR84-b}
N.~Karmarkar.
\newblock A new polynomial-time algorithm for linear programming.
\newblock {\em Combinatorica}, 4(4):373--396, 1984.

\bibitem[Kha79]{KAT79}
L.G. Khachiyan.
\newblock Un algorithme polynomial en programmation linéaire.
\newblock {\em Doklady Akadamii Nauk SSR}, 224:1093--1096, 1979.
\newblock en russe, traduction anglaise: Soviet Math Doklady , Volume 20,
  191-194.

\bibitem[Mur06a]{MUR06-b}
K.G. Murty.
\newblock Linear equations, inequations, lps, and an efficient new algorithm.
\newblock In {\em INFORMS 2006}, pages 3--19, 2006.
\newblock Tutorial.

\bibitem[Mur06b]{MUR06-a}
K.G. Murty.
\newblock A new practically efficient interior point method for lp.
\newblock {\em Algorithmic Operations Research}, 1:3--19, 2006.

\bibitem[NPL09]{NPL09}
N.~C. Nguyen, C.~D. Pham, and T.~H. Le.
\newblock The outer constraction method for linear programming problem.
\newblock Communication personnelle, Hano{\"i}, 2009.

\bibitem[PL09]{PH09}
C.~D. Pham and T.~H. Le.
\newblock An alternating projections algorithm for solving linear programs.
\newblock {\em Acta Mathematica Vietnamica}, 34(3):335--343, 2009.

\bibitem[Tol39]{TOL39}
A.~Tolsto{\"i}.
\newblock M{\'e}thodes d'{\'e}limination des transports non rationnels lors de
  la planification.
\newblock {\em Le transport socialiste}, 1(9):28--51, 1939.
\newblock En Russe.

\end{thebibliography}
%\bibliography{bibio-PL}
%xxxxxxxxxxxxxxxxxxx
%\nocite{*}

\end{document}